\documentclass[12pt]{article}
\usepackage{url}

\usepackage{amsfonts}
\usepackage{amssymb}
\usepackage{acronym}
\usepackage{color}
\usepackage[dvips]{graphicx}
\usepackage{graphicx}
\usepackage{epsfig}

\def\jnt#1{{{#1}}}

\def\ao#1{{{#1}}}
\def\alexo#1{{{#1}}}
\def\an#1{{{#1}}}

\newtheorem{theorem}{Theorem}

\newtheorem{assumption}{Assumption}

\newtheorem{lemma}{Lemma}

\newtheorem{proposition}{Proposition}

\newenvironment{proof}[1][Proof]{\noindent\textbf{#1.} }{\ \rule{0.5em}{0.5em}}
\def\E{\mathcal{E}}
\def\Eb{\overline{\E}}

\def\Gb{\overline{G}}

\def\rn{\mathbb{R}^n}
\def\R{\mathbb{R}}

\begin{document}
\title{ \LARGE \bf On Distributed Averaging Algorithms and Quantization Effects\footnote{This research was partially supported by the National Science Foundation under grants ECCS-0701623, CMMI 07-42538, and DMI-0545910, and by DARPA ITMANET program}}
\author{Angelia Nedi\'c\footnote{A.\ Nedi\'c is with the
Industrial and Enterprise Systems Engineering Department,
University of Illinois at Urbana-Champaign, Urbana IL 61801
(e-mail:angelia@illinois.edu)}, Alex Olshevsky, Asuman Ozdaglar, and John N.\ Tsitsiklis\footnote{
A.\ Olshevsky, A.\ Ozdaglar, and J.\ N.\ Tsitsiklis  are
with the Laboratory for Information and Decision Systems, Electrical Engineering and Computer Science Department, Massachusetts Institute of Technology, Cambridge MA, 02139 (e-mails: alex\_o@mit.edu, asuman@mit.edu, jnt@mit.edu)}}
%\markright{LIDS Report ???}

\maketitle

\thispagestyle{headings}

\begin{abstract}
\noindent We consider distributed iterative algorithms for the
averaging problem over time-varying topologies. Our focus is on the
convergence time of such algorithms when complete (unquantized)
information is available, and on the degradation of performance when
only quantized information is available. We study a large and
natural class of averaging algorithms, which includes the vast
majority of algorithms proposed to date, and provide tight
polynomial bounds on their convergence time. \alexo{We \jnt{also} describe an
algorithm within this class whose convergence time is the best
\jnt{among currently available}
 averaging algorithms for time-varying
topologies.} We then propose and analyze distributed averaging
algorithms under the additional constraint that agents can only
store and communicate quantized information,
\jnt{so that they can only} %. We show that these algorithms
converge to the average of the initial values of the
agents within some error. We establish bounds on the error and tight
bounds on the convergence time, as a function of the number of
quantization levels.
\end{abstract}

\newpage

\section{Introduction}
There has been much recent interest in distributed control and
coordination of networks consisting of multiple, potentially mobile,
agents. This is motivated mainly by the emergence of large scale
networks, characterized by the lack of centralized access to
information and time-varying connectivity. Control and optimization
algorithms deployed in such networks should be completely
distributed, relying only on local observations and information, and
robust against unexpected changes in topology such as link or node
failures.

A canonical problem in distributed control is the {\it consensus
problem}. The objective in the consensus problem is to develop
distributed algorithms that can be used by a group of agents in
order to reach agreement (consensus) on a common decision
(represented by a scalar or a vector value). The agents start with
some different initial decisions and communicate them locally under
some constraints on connectivity and inter-agent information
exchange. The consensus problem arises in a number of applications
including coordination of UAVs (e.g., aligning the agents'
directions of motion), information processing in sensor networks,
and distributed optimization (e.g., agreeing on the estimates of
some unknown parameters). The {\it averaging problem} is a special
case in which the goal is to compute the exact average of the
initial values of the agents. A natural and widely studied consensus
algorithm, proposed and analyzed by Tsitsiklis \cite{T84} and
Tsitsiklis {\it et al.} \cite{TBA86}, involves, at each time step,
every agent taking a weighted average of its own value with values
received from some of the other agents. Similar algorithms have been
studied in the load-balancing literature (see for example
\cite{C89}). Motivated by observed group behavior in biological and
dynamical systems, the recent literature in cooperative control has
studied similar algorithms and proved convergence results under
various assumptions on agent connectivity and information exchange
(see \cite{VCBJCS95}, \cite{OSM04}, \cite{RB05}, \cite{M05},
\cite{LR06}).

In this paper, our goal is to provide tight bounds on the
convergence time (defined as the number of iterations required to
reduce a suitable Lyapunov function by a constant factor) of a
general class of consensus algorithms, as a function of the number
$n$ of agents. We focus on algorithms that are designed to solve the
averaging problem. We consider both problems where agents have
access to exact values and problems where agents only have access
to quantized values of the other agents. Our contributions can be
summarized as follows.

In the first part of the paper, we consider the case where agents
can exchange and store continuous values, which is a widely adopted
assumption in the previous literature. We consider a large class of
averaging algorithms defined by the condition that the weight matrix
is a possibly nonsymmetric, doubly stochastic matrix. For this class
of algorithms, we prove that the convergence time is $O(n^2/\eta)$,
where $n$ is the number of agents and $\eta$ is a lower bound on the
nonzero weights used in the algorithm. To the best of our knowledge,
{\it this is the best polynomial-time bound on the convergence time
of such algorithms}. We also show that this bound is tight. Since
\alexo{all previously studied linear schemes force $\eta$ to be of
the order of $1/n$,} this result implies an $O(n^3)$ bound on
convergence time. In Section \ref{matrixpicking}, we present a
distributed algorithm that selects the weights dynamically, using
three-hop neighborhood information. Under the assumption that the
underlying connectivity graph at each iteration is undirected, we
establish an improved $O(n^2)$ upper bound on convergence time. This
matches the best \ao{currently} available convergence time guarantee
for the much simpler case of static connectivity graphs
\alexo{\cite{OT06}}.

In the second part of the paper, we impose the additional constraint
that agents can only store and transmit quantized values. This model
provides a good approximation for communication networks that are
subject to communication bandwidth \ao{or storage} constraints. We
focus on a particular quantization rule, which rounds down the
values to the nearest quantization level. We propose a distributed
algorithm that uses quantized values and, using a slightly different
Lyapunov function, we show that the algorithm guarantees the
convergence of the values of the agents to a common value. In
particular, we prove that all agents have the same value after
$O((n^2/\eta)\log (nQ))$ time steps, where $Q$ is the number of
quantization levels per unit value. Due to the rounding-down feature
of the quantizer, this algorithm does not preserve the average of
the values at each iteration. However, we provide bounds on the
error between the final consensus value and the initial average, as
a function of the number $Q$ of available quantization levels. In
particular, we show that the error goes to 0 at a rate of $(\log
Q)/Q$, as the number $Q$ of quantization levels increases to
infinity.

Other than the papers cited above, our work is also related to
\cite{KBS06} and \cite{CFSZ05, CFFTZ07}, which studied the effects
of quantization on the performance of averaging algorithms. In
\cite{KBS06}, Kashyap {\it et al.} proposed randomized {\it
gossip-type} quantized averaging algorithms under the assumption
that each agent value is an integer. They showed that these
algorithms preserve the average of the values at each iteration and
converge to approximate consensus. They also provided bounds on the
convergence time of these algorithms for specific static topologies
(fully connected and linear networks). In the recent work
\cite{CFFTZ07}, Carli {\it et al.} proposed a distributed algorithm
that uses quantized values and preserves the average at each
iteration. They showed favorable convergence properties using
simulations on some static topologies, and provided performance
bounds for the limit points of the generated iterates. Our results
on quantized averaging algorithms differ from these works in that
{\it we study a more general case of time-varying topologies, and
provide tight polynomial bounds on both the convergence time and the
discrepancy from the initial average, in terms of the number of
quantization levels}.

The paper is organized as follows. In Section \ref{agreementalg}, we
introduce a general class of averaging algorithms, and present our
assumptions on the algorithm parameters and on the information
exchange among the agents. In Section \ref{convtimesection}, we
present our main result on the convergence time of the averaging
algorithms under consideration. In Section \ref{matrixpicking}, we
present a distributed averaging algorithm for the case of undirected
graphs, which picks the weights dynamically, resulting in an
improved bound on the convergence time. In Section \ref{qanalysis},
we propose and analyze a quantized version of the averaging algorithm.
In particular,
we establish bounds on the convergence time of the iterates, and on the
error between the final value and the average of the initial values
of the agents. Finally, we give our concluding remarks in Section
\ref{conclusions}.

\section{A Class of Averaging Algorithms}\label{agreementalg}

We consider a set $N=\{1,2,\ldots,n \}$ of agents, which will
henceforth be referred to as ``nodes.'' Each node $i$ starts with a
scalar value $x_i(0)$. At each nonnegative integer time $k$, node
$i$ receives from some of the other nodes $j$ a message with the
value of $x_j(k)$,  and updates its value according to:
%the following relation:
\begin{equation}
x_i(k+1) = \sum_{j=1}^n a_{ij}(k) x_j(k), \label{noquant}
\end{equation} where the $a_{ij}(k)$ are nonnegative weights with the
property that $a_{ij}(k)>0$ only if node $i$ receives information
from node $j$ at time $k$. We use the notation $A(k)$ to denote the
{\it weight matrix} $[a_{ij}(k)]_{i,j=1,\ldots,n}$,
\alexo{so that our update equation is \[ x(k+1) = A(k) x(k).\]}
Given a matrix $A$, we use $\E(A)$ to denote the set of directed edges $(j,i)$,
including self-edges $(i,i)$, such that $a_{ij}>0$. At each time
$k$, the nodes' connectivity can be represented by the directed
graph $G(k)=(N,\E(A(k)))$.

Our goal is to study the convergence of the iterates $x_i(k)$ to the
average  of the initial values, $(1/n)\sum_{i=1}^n
x_i(0)$, as $k$ approaches infinity. In order to establish such %a
convergence, we impose some assumptions on the weights $a_{ij}(k)$
and the graph sequence $G(k)$.

\begin{assumption}For each $k$, the weight matrix $A(k)$
 is a doubly stochastic matrix\footnote{\alexo{A matrix is called doubly stochastic
 if it is nonnegative and all of its rows and columns sum to $1$.}} with positive diagonal \alexo{entries}. Additionally, there exists a constant $\eta> 0$ such that if $a_{ij}(k)>0$, then $a_{ij}(k) \geq \eta$.\label{weights}
\end{assumption}

The double stochasticity assumption on the weight matrix guarantees
that the average of the node values remains the same at each
iteration (cf.\ the proof of Lemma \ref{vl} below). The second part
of this assumption states that each node gives significant weight to
its values and to the values of its neighbors at each time $k$.

Our next assumption ensures that the graph sequence $G(k)$ is
sufficiently connected for the nodes to repeatedly influence each
other's values.

\begin{assumption} There exists an integer $B \geq 1$
such that the directed graph \[ \Big(N, \E(A(kB)) \bigcup
\E(A(kB+1))\bigcup \cdots \bigcup \E(A((k+1)B-1))\Big)\]  is strongly
connected for all nonnegative integers $k$.
\label{connectivity}\end{assumption}

Any algorithm of the form given in Eq.\ (\ref{noquant}) with the
sequence of weights $a_{ij}(k)$ satisfying Assumptions \ref{weights}
and \ref{connectivity} solves the averaging problem. This is
formalized in the following \alexo{proposition}.

\begin{proposition} \label{uqc} Let Assumptions \ref{weights} and \ref{connectivity} hold. Let $\{x(k)\}$ be generated by the algorithm (\ref{noquant}). Then,  for all $i$, we have
\[ \lim_{k \rightarrow \infty} x_i(k) = \frac{1}{n} \sum_{j=1}^n
x_j(0). \]
\end{proposition}

This \alexo{fact} is a minor modification of known results in
\cite{T84, TBA86, JLM03, BHOT05}, where the convergence of each
$x_i(k)$ to the same value is established under weaker versions of
Assumptions \ref{weights} and \ref{connectivity}. The fact that the
limit is the average of the entries of the vector $x(0)$ follows
from the fact that multiplication of a vector by a doubly stochastic
matrix preserves the average of the vector's components.

Recent research has focused on methods of choosing weights
$a_{ij}(k)$ that satisfy Assumptions \ref{weights} and
\ref{connectivity}, and minimize the convergence time of the
resulting averaging algorithm (see \cite{XB04} for the case of
static graphs, see \cite{OSM04} and \cite{BFT05} for the case of
symmetric weights, i.e., weights satisfying $a_{ij}(k)=a_{ji}(k)$,
and also see \cite{C06, boyd}). For static graphs, some recent
results on optimal time-invariant algorithms may be found in
\cite{OT06}.

\section{Convergence time \label{convtimesection}}

In this section, we give an analysis of the convergence time of
averaging algorithms of the form (\ref{noquant}). Our goal is to
obtain tight estimates of the convergence time, under Assumptions
\ref{weights} and \ref{connectivity}.

As a convergence measure, we use the ``sample variance'' of a vector
$x \in \R^n$, defined as \[ V(x) = \sum_{i=1}^n (x_i - \bar{x} )^2,
\] where $\bar{x}$ is the average of the entries of $x$: \[ \bar{x} = \frac{1}{n} \sum_{i=1}^n x_i. \]

Let $x(k)$ denote the vector of node values at time $k$ [i.e., the
vector of iterates generated by algorithm (\ref{noquant}) at time
$k$]. We are interested in providing an upper bound on the number of
iterations it takes for the ``sample variance'' $V(x(k))$ to
decrease to a small fraction of its initial value $V(x(0))$. We
first establish some technical preliminaries that will be key in the
subsequent analysis. In particular, in the next subsection, we
explore several implications of the double stochasticity assumption
on the weight matrix $A(k)$.

\subsection{Preliminaries on Doubly Stochastic Matrices}

We begin by analyzing how the sample variance $V(x)$ changes when
the vector $x$ is multiplied by a doubly stochastic matrix $A$. The
next lemma shows that $V(Ax) \leq V(x)$. Thus, under Assumption
\ref{weights}, the sample variance $V(x(k))$ is nonincreasing in
$k$, and $V(x(k))$ can be used as a Lyapunov function.

\begin{lemma}
\label{vl} Let $A$ be a doubly stochastic matrix. Then,\footnote{In
the sequel, the notation $\sum_{i<j}$ will be used to denote the
double sum $\sum_{j=1}^n\sum_{i=1}^{j-1}$.} for all $x \in \R^n$,
\[ V(Ax) = V(x) - \sum_{i<j} w_{ij} (x_i - x_j)^2, \] where $w_{ij}$ is
the $(i,j)$-th entry of the matrix $A^T A$.
\end{lemma}

\begin{proof}
Let ${\bf 1}$ denote the vector in $\R^n$ with all entries equal to
$1$. The double stochasticity of $A$ implies
\[A {\bf 1}= {\bf 1}, \qquad {\bf 1}^T A = {\bf 1}^T. \]
Note that multiplication by a doubly stochastic matrix $A$ preserves
the average of the entries of a vector, i.e., for any $x\in\rn$,
there holds
\[ \overline{Ax} = \frac{1}{n}\, {\bf 1}^T Ax
=\frac{1}{n}\, {\bf 1}^T x = \bar{x}. \] We now write the
quadratic form $V(x)-V(Ax)$ explicitly, as follows:
\begin{eqnarray} V(x) - V(Ax) & = &  (x - \bar{x} {\bf 1})^T (x-\bar{x}{\bf 1}) -
(Ax-\overline{Ax}{\bf 1})^T (Ax - \overline{Ax}{\bf 1}) \nonumber \\
& = &  (x - \bar{x} {\bf 1})^T (x-\bar{x}{\bf 1}) - (Ax-\bar{x}A{\bf
1})^T (Ax - \bar{x} A {\bf 1}) \nonumber \\ & = &  (x - \bar{x} {\bf
1})^T (I-A^T A) (x-\bar{x}{\bf 1}). \label{vdec} \end{eqnarray}

Let $w_{ij}$ be the $(i,j)$-th entry of $A^T A$.  %We next
%parametrize the matrix $A^T A$ in terms of the scalars $w_{ij}$.
Note that $A^T A$ is symmetric and stochastic, so that
$w_{ij}=w_{ji}$ and $w_{ii} = 1 - \sum_{j \neq i} w_{ij}$. Then, it
can be verified that
\begin{equation} A^T A  = I - \sum_{i<j} w_{ij} (e_i - e_j) (e_i -
e_j)^T, \label{asquareform}
\end{equation} where $e_i$ is a unit vector with the $i$-th entry equal to 1,
and all other entries equal to 0 \ao{(see also \cite{XBK07} where a
similar decomposition was used)}.

By combining Eqs.\ (\ref{vdec}) and (\ref{asquareform}), we obtain
\begin{eqnarray*}
V(x) - V(Ax) & =&  (x - \bar{x} {\bf 1})^T
\Big(\sum_{i<j} w_{ij} (e_i - e_j) (e_i - e_j)^T\Big) (x-\bar{x}{\bf 1}) \nonumber \\
& = &  \sum_{i<j} w_{ij} (x_i - x_j)^2. %\label{vdecb}
\end{eqnarray*}
\end{proof}
\vskip 1pc

Note that the entries $w_{ij}(k)$ of $A(k)^T A(k)$ are nonnegative,
because the weight matrix $A(k)$ has nonnegative entries. In view of this,
Lemma \ref{vl} implies that
%Since $x(k+1)=A(k)x(k)$ for each $k$, and the weight matrix $A(k)$
%is doubly stochastic with nonnegative entries [implying that the
%entries $w_{ij}(k)$ of $A(k)^T A(k)$ are nonnegative], it follows
%from the previous lemma that
\[V(x(k+1))\le V(x(k))\qquad \hbox{for all }k.\]
Moreover, the amount of variance decrease is given by
\[V(x(k)) - V(x(k+1)) = \sum_{i<j} w_{ij}(k) (x_i(k) - x_j(k))^2.\]
We will use this result to provide a lower bound on the amount of
decrease of the sample variance $V(x(k))$ in between iterations.

Since every positive entry of $A(k)$ is at least $\eta$, it follows
that every positive entry of $A(k)^T A(k)$ is at least $\eta^2$.
Therefore, it is immediate that
\[ \hbox{if }\quad w_{ij}(k) > 0, \quad\hbox{then }\quad w_{ij}{(k)} \geq \eta^2. \]
In our next lemma, we establish a stronger lower bound. In
particular, we find it useful to focus not on an individual
$w_{ij}$, but rather on all $w_{ij}$ associated with edges $(i,j)$
{that} cross {a particular} cut in the graph $(N,\E(A^TA))$. For
such groups of $w_{ij}$, we prove a lower bound which is linear in
$\eta$, as seen in the following.

\begin{lemma} \label{db}
Let $A$ be a \alexo{row-}stochastic matrix with positive diagonal
\an{entries},
and assume that \an{the smallest positive entry in $A$} is at least $\eta$.
Also, let $(S^{-},S^{+})$ be a partition of the set
$N=\{1,\ldots,n\}$ into two disjoint sets. If
\[ \sum_{i \in S^{-}, ~j \in S^{+}} w_{ij} > 0, \] then
\[ \sum_{i
\in S^{-}, ~j \in S^{+}} w_{ij} \geq \frac{\eta}{2}. \]
\end{lemma}

\begin{proof}
Let $\sum_{i \in S^{-}, ~j \in S^{+}} w_{ij} >0$. From the
definition of the weights $w_{ij}$, we have $w_{ij}=\sum_k
a_{ki}a_{kj}$, which shows that there exist $i\in S^-$, $j\in S^+$,
and some $k$ such that $a_{ki}>0$ and $a_{kj}>0$. For either case
where $k$ belongs to $S^-$ or $S^+$, we see that there exists an
edge in the set $\E(A)$ that crosses the cut $(S^{-},S^{+})$. Let
$(i^*,j^*)$ be such an edge. Without loss of generality, we assume
that $i^* \in S^{-}$ and $j^* \in S^{+}$.
%Consider the edge set $\E(A)$. It must contain an edge crossing the
%cut $(S_{-},S_{+})$; if not, $\E(A^T A)$ would not contain an edge
%crossing the cut $(S_{-},S_{+})$ either, and we would have $\sum_{i
%\in S_{-}, ~j \in S_{+}} w_{ij} = 0$.

%Let $(i^*, j^*)$ be an edge in $\E(A)$ crossing the cut
%$(S_{-},S_{+})$. Without loss of generality, we may assume that $i^*
%\in S_{-}$ and $j^* \in S_{+}$.

We define
\begin{eqnarray*} C_{j^*}^{+} & = & \sum_{i \in S^{+}} a_{j^*i}, \\
C_{j^*}^{-} & =&  \sum_{i \in S^{-}} a_{j^*i}. \end{eqnarray*} See
Figure \ref{cdef}(a) for an illustration.
Since $A$ is a \an{row}-stochastic matrix, we have
\[%\begin{equation}
C_{j^*}^{-} + C_{j^*}^{+}  =  1, %\label{isum}
\]%\end{equation}
implying that at least one of the following is true:
\begin{eqnarray*}
\mbox{ Case (a): } & & C_{j^*}^{-} \geq \frac{1}{2}, \\
\mbox{ Case (b): } & & C_{j^*}^{+} \geq \frac{1}{2}.
\end{eqnarray*}
We consider these {two cases} separately. In both cases,
we focus on a subset of the edges and we use the fact
that the
elements $w_{ij}$ correspond to paths of length $2$, with one step
in $\E(A)$ and another in $\E(A^T)$.

\begin{center}
\begin{figure}
\hspace{1cm}
{\includegraphics[width=4cm]{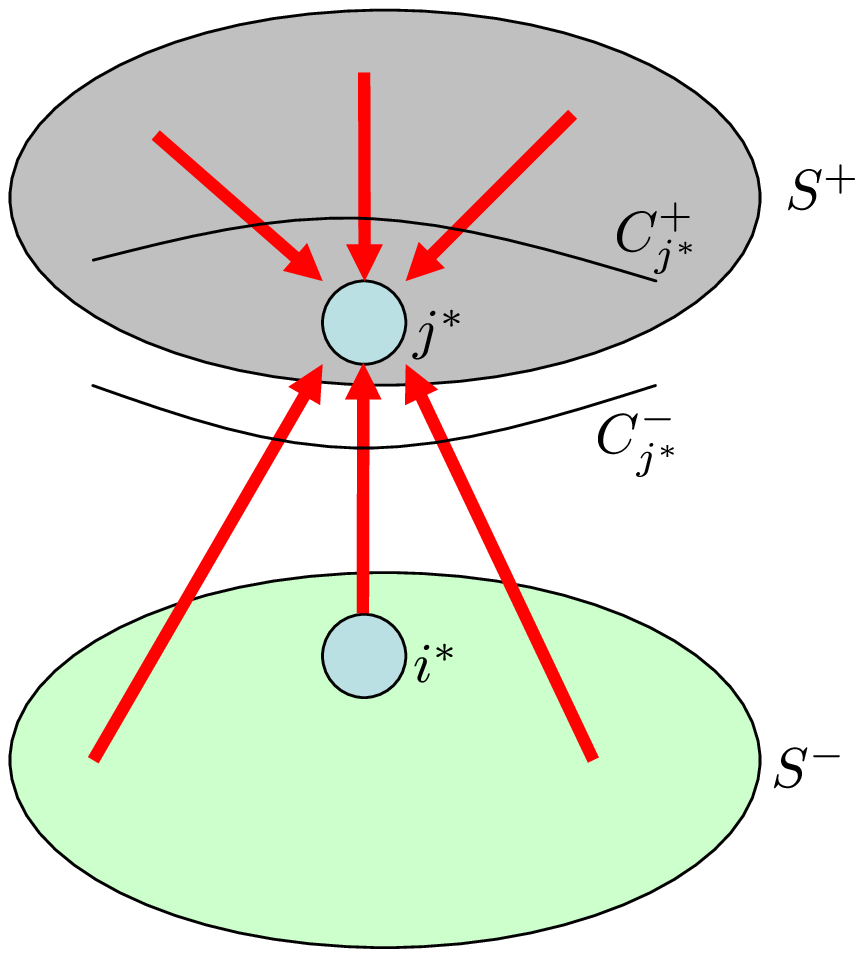} \quad
\includegraphics[width=4cm]{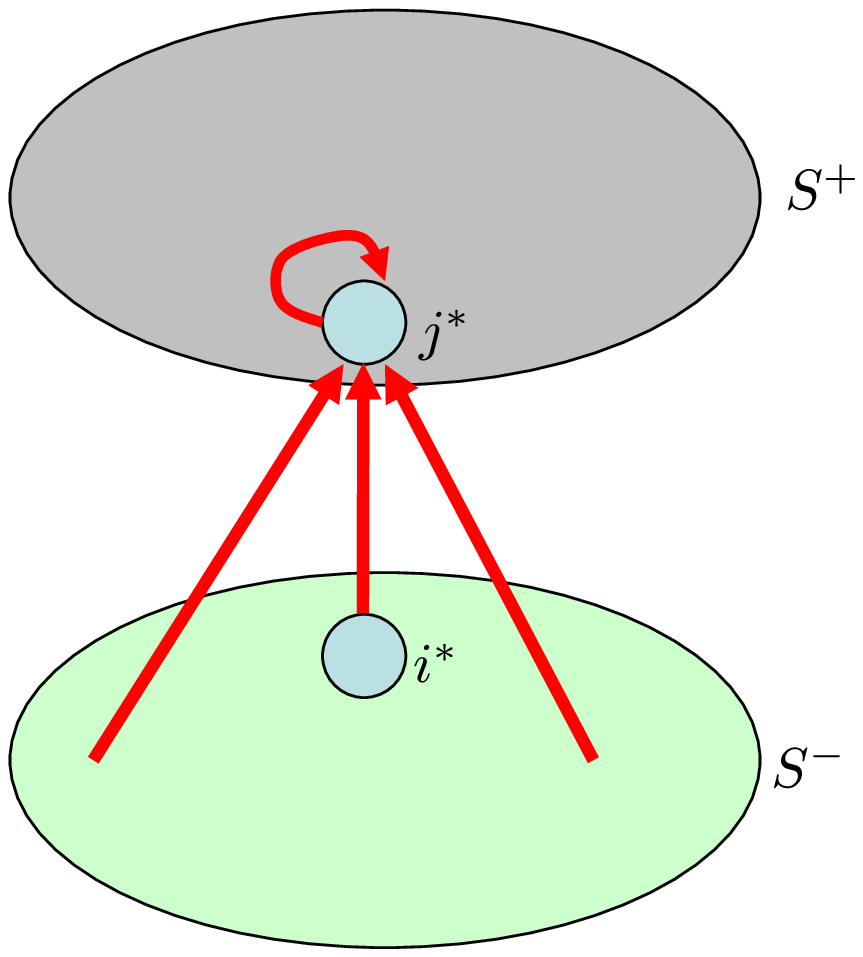}
\quad
\includegraphics[width=4cm]{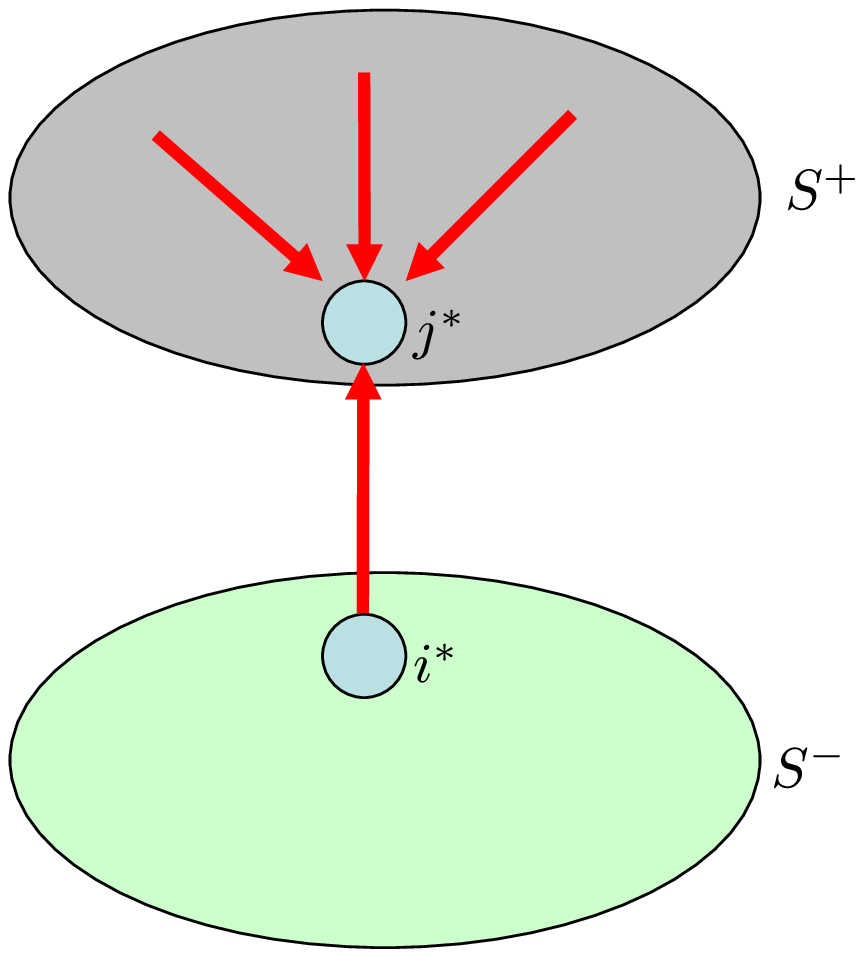}
} \caption{\label{cdef}\small (a) Intuitively,  $C_{j^*}^{+}$
measures how much weight $j^*$ assigns to nodes in $S^{+}$ (including
itself), and $C_{j^*}^{-}$ measures how much weight $j^*$ assigns to
nodes in $S^{-}$. Note that the {edge} $(j^*,j^*)$ is also present,
but not shown.
(b) For the case {where} $C_{j^*}^{-} \geq 1/2$, we
only focus on
two-hop paths between $j^*$ and elements $i \in S^{-}$
obtained by taking $(i,j^*)$ as the first step and the self-edge
$(j^*,j^*)$ as the second step.
(c) For the case where
$C_{j^*}^{+} \geq 1/2$, we
only focus on two-hop paths between $i^*$ and elements $j\in S^{+}$
obtained by taking $(i^*,j^*)$ as the first step in
$\E(A)$ and $(j^*,j)$ as the second step in $\E(A^T)$.}
\end{figure}
\end{center}

\noindent {\it Case (a):} $C_{j^*}^{-} \geq 1/2.$\\
We focus on those $w_{ij}$ with $i\in S^-$ and $j=j^*$. Indeed,
since all $w_{ij}$ are nonnegative, we have
\begin{equation} \label{wsubset} \sum_{i \in S^{-}, ~j \in S^{+}}
w_{ij} \geq \sum_{i \in S^{-}} w_{ij^*}.
\end{equation}
For each element in the sum on the right-hand {side}, we have
\[ w_{ij^*} = \sum_{k=1}^n a_{ki}\ a_{kj^*}
\geq a_{j^*i}\ a_{j^* j^*} \geq a_{j^*i}\ \eta, \] where the
inequalities follow from the facts that $A$ has
nonnegative entries, its diagonal entries are positive, and
its positive entries are at least $\eta$. Consequently,
\begin{equation} \sum_{i \in S^{-}} w_{ij^*}
\geq \eta\ \sum_{i \in S^{-}} a_{j^*i}
= \eta\, C_{j^*}^{-}. \label{wlowerbound}
\end{equation}
Combining Eqs.\ (\ref{wsubset}) and (\ref{wlowerbound}), and
recalling {the assumption} $C_{j^*}^{-} \geq 1/2$, the {result}
follows. {An} illustration of {this argument} can be found in Figure
\ref{cdef}{(b)}.

\noindent {\it Case (b):} $C_{j^*}^{+} \geq 1/2.$\\
We focus on those $w_{ij}$ with $i=i^*$ and $j \in S^{+}$.
We have
\begin{equation} \label{wsubset2}
\sum_{i \in S_{-}, ~j \in S^{+}} w_{ij} \geq \sum_{j \in
S^{+}} w_{i^* j},
\end{equation}
since all $w_{ij}$ are nonnegative. For each element in the sum on
the right-hand side, we have
\[ w_{i^* j} = \sum_{k=1}^n a_{ki^*}\ a_{k j}
\geq a_{j^* i^*}\ a_{j^* j} \geq \eta\, a_{j^* j}, \]
where the inequalities follow {because all entries of $A$ are
nonnegative, and because the choice $(i^*, j^*) \in \E(A)$ implies
that $a_{j^* i^*} \geq{\eta}$.} Consequently,
\begin{equation}
\sum_{j \in S^{+}} w_{i^* j} \geq \eta \sum_{j \in
S^{+}} a_{j^*j} = \eta\, C_{j^*}^{+}. \label{wlowerbound2}
\end{equation} {Combining Eqs.\ (\ref{wsubset2}) and}
(\ref{wlowerbound2}), and recalling {the assumption} $C_{j^*}^+ \geq
1/2$, the {result} follows. An illustration of this argument can
be found in Figure \ref{cdef}{(c)}.
\end{proof}

\subsection{A Bound on Convergence Time}\label{convergtime}

With the preliminaries on doubly stochastic matrices in place, we
{can now proceed to derive} bounds on the decrease of $V(x(k))$ in
between iterations. We will first somewhat relax our connectivity
assumptions. In particular, we consider the following relaxation of
Assumption \ref{connectivity}.\vspace{5pt}

\begin{assumption} \label{weakconnect} Given an integer $k\ge 0$,
suppose that the components of $x(kB)$ have been reordered so that
they are in nonincreasing order. We assume that for every
$d\in\{1,\ldots,n-1\}$, we either have $x_d(kB)=x_{d+1}(kB)$, or
there exist some time $t\in\{kB,\ldots,(k+1)B-1\}$ and some
$i\in\{1,\ldots,d\}$, $j\in\{d+1,\ldots,n\}$ such that $(i,j)$ or
$(j,i)$ belongs to $\E(A(t))$.
\end{assumption}

\begin{lemma} \label{as2implas3}
Assumption \ref{connectivity}
implies Assumption \ref{weakconnect}, {with the same value of $B$.}
\end{lemma}

\begin{proof} If Assumption \ref{weakconnect} does not hold,
then \ao{there must exist an index $d$
[for which $x_d(kB) \neq x_{d+1}(kB)$ holds] such that there are no
edges between nodes $1,2,\ldots,d$ and nodes $d+1,\ldots,n$ during
times $t=kB,\ldots,(k+1)B-1$. But this implies that the graph}
\[ \Big(N,\E({A}(kB)) \bigcup \E(A(kB+1)){\bigcup} \cdots \bigcup
\E(A((k+1)B-1))\Big) \]
is disconnected, which violates Assumption 2.
\end{proof}

For our convergence time results, we will use the weaker Assumption
\ref{weakconnect}, rather than the stronger Assumption 2. Later on,
in Section \ref{matrixpicking}, we will exploit the sufficiency of
Assumption \ref{weakconnect} to {design} a % propose a constructive,
decentralized algorithm for selecting the weights $a_{ij}(k)$, which
satisfies Assumption \ref{weakconnect}, but not  Assumption 2.

We now proceed to bound the decrease of our Lyapunov function
$V(x(k))$ during the interval $[kB, (k+1) B-1]$. In
what follows, we denote by $V(k)$ the sample variance $V(x(k))$ at
time $k$.

\begin{lemma}
\label{vardiff} Let Assumptions \ref{weights} and \ref{weakconnect}
hold. Let $\{x(k)\}$ be generated by the update rule
(\ref{noquant}).
Suppose that the components $x_i(kB)$ of the vector
$x(kB)$
have been ordered from largest to smallest, with ties broken
arbitrarily. Then, we have
\[%\begin{equation} \label{vardiffeq}
V(kB) - V((k+1)B) \geq
\frac{\eta}{2} \sum_{i=1}^{n-1} (x_{i}(kB) - x_{i+1}(kB))^2.
\]%\end{equation}
\end{lemma}

\begin{proof} By Lemma \ref{vl}, we have for all $t$,
\begin{equation}
V(t)-V(t+1) = \sum_{i<j} w_{ij}(t)
(x_i(t)-x_j(t))^2, \label{firstvdec0}
\end{equation}
where $w_{ij}(t)$ is the $(i,j)$-th entry of $A(t)^T A(t)$. Summing
up the variance differences $V(t)-V(t+1)$ over different values of $t$, we obtain
\begin{equation}
V(kB) - V((k+1)B)
= \sum_{t=kB}^{(k+1)B-1} \sum_{i<j}
w_{ij}(t) (x_i(t)-x_j(t))^2. \label{firstvdec}
\end{equation}

\ao{We next introduce some notation.}

\begin{itemize}
\item[(a)]
For all $d\in\{1,\ldots,n-1\}$, let $t_d$ be the first time larger than
or equal to $kB$
(if it exists) at which there is a communication between two nodes
belonging to the two sets $\{1,\ldots,d\}$ and $\{d+1,\ldots,n\}$, to
be referred to as a communication across the cut $d$.

\item[(b)]
For all $t\in \{kB,\ldots,(k+1)B-1\}$, let $D(t)=\{d\mid
t_d=t\}$, i.e., $D(t)$ consists of ``cuts" $d\in \{1,\ldots,n-1\}$
such that time $t$ is the first communication time larger than or equal
to $kB$ between nodes in the sets $\{1,\ldots,d\}$ and
$\{d+1,\ldots,n\}$. Because of Assumption \ref{weakconnect}, the
union of the sets $D(t)$ includes all indices $1,\ldots,n-1$, except
possibly for indices for which $x_d(kB)=x_{d+1}(kB)$.

\item[(c)]
For all $d\in \{1,\ldots,n-1\}$, let $C_d=\{(i,j) \alexo{,~(j,i)} %\ao{\in
%\E(A(t_d))}
 \mid i\leq d, \ d+1\leq j\}$.

\item[(d)]
For all $t\in \{kB,\ldots,(k+1)B-1\}$, let $F_{ij}(t) =\{d\in D(t)\
|\ (i,j) \alexo{\mbox{ or } (j,i)} \in C_d\}$, i.e., $F_{ij}(t)$
consists of all cuts $d$ such that the edge $(i,j)$ \alexo{ or
$(j,i)$} at time $t$ is the first communication across the cut at a
time larger than or equal to $kB$.

\item[(e)] To simplify notation, let $y_i=x_i(kB)$. By assumption, we have
$y_1\geq\cdots\geq y_n$.
\end{itemize}

We make two observations, as follows:
\begin{itemize}
\item[(1)]
Suppose that $d\in D(t)$. Then, for some $(i,j)\in C_d$, we have
either $a_{ij}(t)>0$ or $a_{ji}(t)>0$.
\alexo{Because $A(t)$ is
nonnegative with positive diagonal \an{entries, we have}
\[ w_{ij}(t) = \sum_{k=1}^n
a_{ki} a_{kj} \geq a_{ii}(t) a_{ij}(t) + a_{ji}(t) a_{jj}(t) > 0,\]
and} by Lemma \ref{db}, we obtain
\begin{equation}
\sum_{(i,j)\in C_d} w_{ij}(t)\geq \frac{\eta}{2}. \label{eq:w}
\end{equation}
\item[(2)]
Fix some $(i,j)$, with $i<j$, and time $t\in
\{kB,\ldots,(k+1)B-1\}$, and suppose that $F_{ij}(t)$ is nonempty.
Let $F_{ij}\ao{(t)}=\{d_1,\ldots,d_k\}$, where the $d_j$ are
arranged in increasing order. Since $d_1\in F_{ij}\ao{(t)}$, we have
$d_1\in D(t)$ and therefore $t_{d_1}=t$. By the definition of
$t_{d_1}$, this implies that there has been no communication between
a node in $\{1,\ldots,d_1\}$ and a node in $\{d_1+1,\ldots,n\}$
during the time interval $[kB,t-1]$. It follows that $x_i(t) \geq
y_{d_1}$. By a symmetrical argument, we also have
\begin{equation}
x_j(t)\leq y_{d_k+1}.\label{cutbound}
\end{equation}
These relations imply that
$$x_i(t)-x_j(t) \geq y_{d_1}-y_{d_k+1}  \ao{\geq } \sum_{d\in F_{ij}\ao{(t)}} (y_d-y_{d+1}),$$
 Since the components of $y$ are sorted in nonincreasing order, we
have $y_d-y_{d+1}\geq 0$, for every $d\in F_{ij}\ao{(t)}$. For any
nonnegative numbers $z_i$, we have
%Using the fact that when the numbers $z_i$ are nonnegative, we have
$$(z_1+\cdots+z_k)^2\geq z_1^2+\cdots+z_k^2,$$
which implies that %and we conclude
\begin{equation}
(x_i(t)-x_j(t))^2 \geq \sum_{d\in F_{ij}\ao{(t)}} (y_d-y_{d+1})^2.
\label{eq:dx}
\end{equation}
%\item[(3)] \alexo{Similarly, one can argue that if $F_{ji}(t)$ is nonempty,
%then \begin{equation} (x_i(t)-x_j(t))^2 \geq \sum_{d\in
%F_{ji}\ao{(t)}} (y_d-y_{d+1})^2. \label{eq:dx2}
%\end{equation}}
\end{itemize}
We now use these two observations to provide a lower bound on the
expression on the right-hand side of Eq.\ (\ref{firstvdec0}) at time
$t$.
%by focusing on the edges that go from nodes in $\{1,\ldots,d\}$
%to nodes in $\{d+1,\ldots,n\}$ for cuts $d$ with $t_d=t$. More
%specifically,
We use Eq.\ (\ref{eq:dx}) and then Eq.\ (\ref{eq:w}),
to obtain
\begin{eqnarray*}
\sum_{i<j} w_{ij}(t)(x_i(t)-x_j(t))^2&
\geq& \sum_{i<j} w_{ij}(t) \sum_{d\in F_{ij}\ao{(t)}} (y_d-y_{d+1})^2\\
&=&\sum_{d\in D(t)} \sum_{(i,j)\in C_d} w_{ij}(t) (y_d-y_{d+1})^2\\
&\geq& \frac{\eta}{2} \sum_{d\in D(t)}(y_d-y_{d+1})^2.
\end{eqnarray*}
We now sum both sides of the above inequality for different values of $t$, and
use Eq.\ (\ref{firstvdec}), to obtain
\begin{eqnarray*}
V(kB)-V((k+1)B)&=& \sum_{t=kB}^{(k+1)B-1} \sum_{i<j} w_{ij}(t)(x_i(t)-x_j(t))^2\\
&\geq&
\frac{\eta}{2} \sum_{t=kB}^{(k+1)B-1}\sum_{d\in D(t)}(y_d-y_{d+1})^2\\
&=&\frac{\eta}{2}\sum_{d=1}^{n-1} (y_d-y_{d+1})^2,
\end{eqnarray*}
where the last inequality follows from the fact that the union of
the sets $D(t)$ is only missing those $d$ for which $y_d=y_{d+1}$.
\end{proof}

\vskip 1pc We next establish a bound on {the} variance decrease that
plays {a} key role in our convergence analysis.

\begin{lemma} \label{lboundvar} Let Assumptions
\ref{weights} and \ref{weakconnect} hold, and suppose that
$V(kB)>0$. Then,
\[\frac{V(kB)-V((k+1)B)}{V(kB)} \geq \frac{\eta}{2 n^2}\qquad
\hbox{for all }k.\]
\end{lemma}

\begin{proof}
{Without loss of generality, we assume that the components of
$x(kB)$ have been sorted in nonincreasing order.} By Lemma
\ref{vardiff}, we have
\[ V(kB) -
V((k+1)B) \geq \frac{\eta}{2}\, \sum_{i=1}^{n-1} (x_{i}(kB) -
x_{i+1}(kB))^2.\]
This implies that
\[ \frac{V(kB)-V((k+1)B)}{V(kB)}\ge\frac{\eta}{2} \,
\frac{\sum_{i=1}^{n-1} (x_{i}(kB) - x_{i+1}(kB))^2}{\sum_{i=1}^n
(x_i(kB)-\bar{x}(kB))^2}. \] Observe that the right-hand side does
not change when we add a constant to every $x_i(kB)$. We can
therefore assume, without loss of generality, that $\bar{x}(kB)=0$,
{so that}
\[ \frac{V(kB)-V((k+1)B)}{V(kB)} \geq \frac{\eta}{2}\,
\min_{{x_1\ge x_2\ge \cdots\ge x_n \atop \sum_i x_i=0}}
\frac{\sum_{i=1}^{n-1} (x_{i} - x_{i+1})^2}{\sum_{i=1}^n x_i^2}. \]
Note that the right-hand side is unchanged if we multiply each
$x_i$ by the same constant. Therefore, we can assume, without loss
of generality, that $\sum_{i=1}^n x_i^2=1$, {so that}
\begin{equation}\label{eq:vv}
\frac{V(kB)-V((k+1)B)}{V(kB)} \geq \frac{\eta}{2}\, \min_{{x_1\ge
x_2\ge \cdots\ge x_n  \atop \sum_i x_i=0, \ \sum_i x_i^2 = 1}}
~~\sum_{i=1}^{n-1} (x_{i} - x_{i+1})^2.\end{equation} The
requirement $\sum_i x_i^2 = 1$ implies that the average value of
$x_i^2$ is $1/n$, which implies that there exists some $j$ such
that $|x_j| \ge 1/\sqrt{n}$. Without loss of generality, let us
suppose that this $x_j$ is positive.\footnote{Otherwise, we can
replace $x$ with $-x$ and subsequently reorder to maintain the
property that the components of $x$ are in descending order. It can
be seen that these operations do not affect the objective value.}

The rest of the proof relies on a technique from \cite{LO81} to
provide a lower bound on the right-hand side of Eq.\ (\ref{eq:vv}).
Let
\[z_i = x_{i} - x_{i+1} \ \ \hbox{for } i<n,\quad
\hbox{and}\quad z_n=0.\] Note that $z_i \geq 0$ for all $i$ \ao{and
\[ \sum_{i=1}^n z_i = x_1 - x_n.\]}Since $x_j\ge 1/\sqrt{n}$ for
some $j$, \ao{we have that $x_1 \geq 1/\sqrt{n}$; since
$\sum_{i=1}^n x_i = 0$, it follows that at least one $x_i$ is
negative, and therefore $x_n < 0$. This gives us \[ \sum_{i=1}^n z_i
\geq \frac{1}{\sqrt{n}} .\]} Combining with {Eq.\ (\ref{eq:vv})}, we
obtain
\[ \frac{V(kB)-V((k+1)B)}{V(kB)} \geq
\frac{\eta}{2} \min_{z_i \geq 0,\ \sum_i z_i \geq 1/\sqrt{n}}
\sum_{i=1}^{n} z_i^2. \] {The minimization problem on the right-hand
side is a symmetric convex optimization problem, and therefore has a
symmetric optimal solution, namely} $z_i = 1/n^{1.5}$ for all $i$.
This results in an optimal value of $1/n^2$. Therefore,
\[\frac{V(kB)-V((k+1)B)}{V(kB)} \geq \frac{\eta}{2 n^2}, \]
which {is} the desired result. \end{proof}

\vskip 1pc

We are now ready for our main result, which establishes that the
convergence time of the sequence of vectors $x(k)$ generated by Eq.\
(\ref{noquant}) is of order $O(n^2B/\eta)$.

\noindent \begin{theorem} \label{uqbound} Let Assumptions
\ref{weights} and \ref{weakconnect} hold. Then, there exists an absolute
constant\footnote{\ao{We say $c$ is an absolute constant when it
does not depend on any of the parameters in the problem, in this
case $n,B, \eta,\epsilon$.}} $c$ such that we have
\[V({k})
\leq \epsilon V(0)\qquad \hbox{for all } k\ge {c} (n^2/\eta)B \log
(1/\epsilon).\]
\end{theorem}

\begin{proof}
The result follows immediately from Lemma \ref{lboundvar}.
\end{proof}

\vskip 1pc
Recall that, according to  Lemma \ref{as2implas3},
Assumption \ref{connectivity} implies Assumption
\ref{weakconnect}.
In view of this, the convergence time bound of Theorem \ref{uqbound}
holds for any $n$ and any sequence of weights satisfying
Assumptions \ref{weights} and \ref{connectivity}.
In the next subsection, we show that this bound is tight when the
stronger Assumption \ref{connectivity} holds.

\subsection{Tightness}

The next \alexo{proposition} shows that the convergence time bound
of Theorem \ref{uqbound} is tight under Assumption
\ref{connectivity}.

\begin{proposition} \label{unquanttight}
{There exist constants $c$ and $n_0$ with the following property.
For any $n\geq n_0$,} nonnegative integer $B$, $\eta < 1/2$, and
$\epsilon < 1$, there exist a sequence of weight matrices $A(k)$ satisfying
Assumptions \ref{weights} and \ref{connectivity}, and an initial
value $x(0)$ such that {if} $V(k)/V(0) \leq \epsilon$, {then}
\[ k \geq {c}\, \frac{n^2}{\eta}\, B \log \frac{1}{\epsilon}.
\]
\end{proposition}

\begin{proof} Let $P$ be the circulant shift operator
defined by $P e_i = e_{i+1}$, $P e_n = e_1$, where $e_i$ is a unit
vector with the $i$-th entry equal to 1, and all other entries equal
to 0. Consider the symmetric circulant matrix defined by
\[ A = (1 - 2 \eta) I + \eta P + \eta P^{-1}. \]
Let $A({k})=A$, when ${k}$ is a multiple of $B$, and $A({k})=I$
otherwise. Note that this sequence satisfies Assumptions 1 and 2.

The second largest eigenvalue of $A$ is
\[ \lambda_2(A) = 1 - 2 \eta + 2 \eta \cos \frac{2 \pi}{n},\]
(see Eq.\ (3.7) of
\cite{G06}). Therefore, using the inequality $\cos x \geq 1- x^2/2$,
\[ \lambda_2(A) \geq 1 - \frac{4 \eta \pi^2}{n^2}. \] {For $n$ large enough,} the quantity
on the right-hand side is nonnegative. Let the initial {vector}
$x(0)$ be the eigenvector corresponding to {$\lambda_2(A)$.} Then,
\[ \frac{V(kB)}{V(0)} = \lambda_2(A)^{2k} \geq \Big(1 - \frac{8 \eta
\pi^2}{n^2}\Big)^k. \] {For the right-hand side to become} less than
$\epsilon$, we need  $k = \Omega( (n^2/\eta) \log (1/\epsilon))$.
{This implies that for $V(k)/V(0)$ to become less than $\epsilon$,
we need} $k=\Omega( (n^2/\eta) B \log (1/\epsilon))$.
\end{proof}

\section{Saving a {factor of} $n$: faster averaging on undirected graphs}
\label{matrixpicking}

In the previous section, we have shown that a large class of
averaging algorithms have  $O(B (n^2/\eta) ~\alexo{ \log
1/\epsilon}) $ \alexo{convergence time. Moreover, we have shown that
this bound is tight, in the sense that there exist matrices
satisfying Assumptions \ref{weights} and \ref{weakconnect} which
converge in $\Omega(B (n^2 /\eta) ~\alexo{ \log 1/\epsilon}) $.}

In this section, we consider decentralized ways of synthesizing the
weights $a_{ij}(k)$ while satisfying Assumptions \ref{weights} and
\ref{weakconnect}. \alexo{Our focus is on improving convergence time
bounds by \an{constructing} ``good'' schemes.}

We assume that the communications of the nodes are governed by an
exogenous sequence of graphs $\Gb(k)=(N,\Eb(k))$ that provides
strong connectivity over time periods of length $B$. \alexo{This
sequence of graphs \jnt{constrains} the matrices $A(k)$ \jnt{that we can}
%which we may
use;} in particular, we require that $a_{ij}(k)=0$ if
$(j,i)\notin\Eb(k)$. Naturally, we assume that $(i,i)\in\Eb(k)$ for
every $i$.

Several such decentralized protocols exist. For example,  each node
may assign
\begin{eqnarray*} a_{ij}{(k)} & = & \epsilon, \qquad\qquad \mbox{ \ \ if
} (j,i)
\in \Eb(k) {\mbox{ and } i\neq j,} \\
a_{ii}{(k)} & = & 1 - \epsilon\cdot \mbox{deg}(i),
\end{eqnarray*}
where {\rm deg}($i$) is the degree of $i$ in $\Gb(k)$. If $\epsilon$
is small enough and the graph $\Gb(k)$ is undirected {[i.e.,
$(i,j)\in\Eb(k)$ if and only if $(j,i)\in\Eb(k)$],} this results in
a nonnegative, doubly stochastic matrix (see \cite{OSM04}). However,
{if a node has $\ao{\Theta}(n)$ neighbors, $\eta$ will be of order
$\ao{\Theta}(1/n)$, resulting in $\ao{\Theta}(n^3)$ convergence
time.} Moreover, this
argument applies to all protocols in which nodes assign equal
weights to all their neighbors; see \cite{XB04} and \cite{BFT05} for
more examples.

In this section, we examine whether it is possible to synthesize the
weights $a_{ij}(k)$ in a decentralized manner, so that
$a_{ij}(k)\geq \eta$ whenever $a_{ij}(k)\neq 0$,
where $\eta$ is a positive constant independent of $n$ and $B$.
We show that this is indeed possible, under the additional assumption that
the graphs $\Gb(k)$ are undirected. Our algorithm is data-dependent,
in that $a_{ij}(k)$ depends not only on the graph $\Gb(k)$, but
also on the data vector $x(k)$. Furthermore,
it is a decentralized 3-hop algorithm, in that $a_{ij}{(k)}$ depends only
on the data at nodes within a distance of at most $3$ from $i$.
Our algorithm is such that the
resulting sequences of vectors $x(k)$ and
graphs $G(k)=(N,\E(k))$, with $\E(k)=\{(j,i)\mid a_{ij}(k)> 0\}$,
satisfy Assumptions
\ref{weights} and \ref{weakconnect}. Thus, a convergence time
result can be obtained from Theorem \ref{uqbound}.

\subsection{The algorithm}

The algorithm we present here is a variation of an old {\em load
balancing} algorithm (see \cite{C89} and Chapter 7.3 of
\cite{BT89}).\footnote{This algorithm was also considered in
\cite{OT06}, but in the absence of a result such as Theorem
\ref{uqbound}, a weaker convergence time bound was derived.}

At each step of the algorithm, each node offers some of its value to
its neighbors, and accepts or rejects such offers from its
neighbors. Once an offer from $i$ to $j$, of size $\delta>0$, has
been accepted, the updates $x_i \leftarrow x_i - \delta$ and $x_j
\leftarrow x_j + \delta$ are executed.

We next describe the formal steps the nodes execute at each time
${k}$. For clarity, we refer to the node executing the steps below
as node $\ao{C}$. Moreover, the instructions below sometimes refer
to the neighbors of node $\ao{C}$; this always means current
neighbors at time $k$, when the step is being executed, {as
determined by the current graph $\Gb(k)$.}
% (since $G(t)$ is not
%assumed to be constant with $t$,  the set of neighbors of $A$ may be
%time-varying).
We assume that at each time $k$, all nodes execute these steps in
the order described below, while the graph remains unchanged.

\vspace{1pc}
\noindent {\bf Balancing Algorithm:}
\begin{enumerate}\itemsep=0pt
\item[1.]  Node $\ao{C}$ broadcasts its current value $x_{\ao{C}}$ to all its
neighbors.

\item[2.]  Going through the values it just received from its
neighbors, Node $\ao{C}$ finds the smallest value that is less than
{its own}. Let $\ao{D}$ be a neighbor with this value. Node $\ao{C}$
makes an offer of $(x_{\ao{C}} - x_{\ao{D}})/3$ to node $D$.

If no {neighbor of $\ao{C}$} has a value smaller than $x_{\ao{C}}$,
node $C$ does nothing at this stage.

\item[3.] Node $\ao{C}$ goes through the incoming offers. It sends an
acceptance to the sender of {a} largest offer, and a rejection to
all the other senders. It updates the value of $x_{\ao{C}}$ by
adding the value of the accepted offer.

If node $C$ did not receive any offers, it does nothing at this
stage.

\item[4.] If an acceptance arrives {for} the offer made by node
$\ao{C}$, node $\ao{C}$ updates $x_{\ao{C}}$ by subtracting the
value of the offer.
\end{enumerate}

Note that the new value of each node is a linear combination of the
values of its neighbors. Furthermore, the weights $a_{ij}(k)$ are
completely determined by the data and the graph at most $3$ hops
from node $i$ in $\Gb(k)$. \alexo{This is true because in the course
of execution of the above steps, each node makes at most three
transmission to its neighbors, so the new value of node $C$ cannot
depend on information more than $3$ hops away from $C$.}

%Indeed, note that the offers received by a node $\ao{D}$ depend on
%$\ao{D}$'s neighbors and on its neighbors' neighbors, i.e., on the
%data within $2$ hops. Thus, whether a node $\ao{D}$ will accept an
%offer from node $\ao{C}$ depends on data at most $2$ hops away from
%$\ao{D}$, which is at most $3$ hops away from $\ao{C}$.

\subsection{Performance analysis}

\alexo{In the following theorem, we are able to
\an{remove} a factor of
$n$ \jnt{from} %relative to
the worst-case convergence time bounds of Theorem
\ref{uqbound}}.

\noindent \begin{theorem} \label{savingn}
Consider the balancing algorithm,
and suppose that $\Gb(k)=(N,\Eb(k))$ is a sequence of
undirected graphs such that $(N,\Eb(kB)\cup\Eb(kB+1)\cup \cdots \cup
\Eb((k+1)B-1))$ is connected, for all integers $k$. There exists an
absolute constant $c$ such that we have
\[V({k}) \leq \epsilon V(0)\qquad \hbox{for all }k\ge {c}
n^2B \log (1/\epsilon).\]
\end{theorem}

\begin{proof}
Note that with this algorithm, the new value at some node $i$ is a
convex combination of the previous values of itself and its
neighbors. Furthermore, the algorithm keeps the sum of the nodes'
values constant, because every accepted offer involves an increase
at the receiving node equal to the decrease at the offering node.
These two properties imply that the algorithm can be written in the
form
\[ x(k+1) = A(k) x(k), \] where $A(k)$ is a doubly stochastic matrix,
determined by $\Gb(k)$ and $x(k)$. It can be seen that the diagonal
entries of $A(k)$ are positive and, furthermore, all nonzero entries
of $A(k)$ are larger than or equal to 1/3; thus, $\eta=1/3$.

We claim that the algorithm [in particular, the sequence $\E(A(k))$]
satisfies Assumption \ref{weakconnect}. Indeed, suppose that at time
$kB$, the nodes are reordered so that the values $x_i(kB)$ are
nonincreasing in $i$. Fix some $d\in\{1,\ldots,n-1\}$, and suppose
that $x_d(kB)\neq x_{d+1}(kB)$. Let $S^+=\{1,\ldots,d\}$ and
$S^-=\{d+1,\ldots,n\}$.

Because of our assumptions on the graphs $\Gb(k)$, there will be a
first time $t$ in the interval $\{kB,\ldots,(k+1)B-1\}$, at which
there is an edge in $\Eb(t)$ between some $i^*\in S^+$  and $j^*\in
S^-$. Note that between times $kB$ and $t$, the two sets of nodes,
$S^+$ and $S^-$, do not interact, which implies that $x_i(t)\geq
x_d(kB)$, for $i\in S^+$, and $x_j(t) < x_d(kB)$, for $j\in S^-$.

At time $t$, node $i^*$ sends an offer to a neighbor with the
smallest value; let us denote that neighbor by $k^*$. Since
$(i^*,j^*)\in \Eb(t)$, we have $x_{k^*}(t)\leq x_{j^*}(t)<x_d(kB)$,
which implies that $k^*\in S^-$. Node $k^*$ will accept the largest
offer it receives, which must come from a node with a value no
smaller than $x_{i^*}(t)$, \ao{and therefore no smaller than
$x_d(kB)$}; hence the latter node belongs to $S^+$. It follows that
$\E(A(t))$ contains an edge between $k^*$ and some node in $S^{+}$,
showing that Assumption \ref{weakconnect} is satisfied.

The claimed result follows from Theorem \ref{uqbound}, because we
have shown that all of the assumptions in that theorem are satisfied
\alexo{with $\eta=1/3$}.
\end{proof}

\section{Quantization Effects \label{qanalysis}}

In this section, we consider a quantized version of the update rule
(\ref{noquant}). This model is a good approximation for a network of
nodes communicating through {finite bandwidth channels, so that} at
each time instant, only a finite number of bits can be transmitted.
We incorporate this constraint in our algorithm by assuming that
each node, upon receiving the values of its neighbors, computes the
convex combination $\sum_{j=1}^{{n}} a_{ij}(k) x_j(k)$ and quantizes
it. This update rule also captures {a} constraint that each node
can only store quantized values.

Unfortunately, {under Assumptions \ref{weights} and
\ref{connectivity}, if the output of Eq.\ (\ref{noquant}) is rounded
to the nearest integer, the sequence $x(k)$} is not guaranteed to
converge to consensus; see \cite{KBS06}. We therefore choose a
quantization rule that rounds
the values down, %consistently
according to
\begin{equation} \label{quantupdate}
x_i(k+1) = \left\lfloor \sum_{j=1}^{{n}} a_{ij}(k) x_j(k)
\right\rfloor,
\end{equation}
where $\lfloor \cdot \rfloor$ represents rounding {\it down} to the
nearest multiple of $1/Q$, and where $Q$ is some positive
integer.

We adopt the natural assumption that the initial values are
{already} quantized.

\begin{assumption}
For all $i$, $x_i(0)$ is a multiple of
$1/Q$.\label{quantizedinitials}
\end{assumption}

For convenience we define \[ U = \max_i x_i(0), \qquad L = \min_i
x_i(0). \] We use $K$ to denote the total number of {relevant}
quantization levels, i.e., \[ K = (U-L)Q, \] which is an integer by
Assumption \ref{quantizedinitials}.

\subsection{A quantization level dependent bound}

We first present a {convergence time bound that depends on the
quantization level $Q$.}

\begin{proposition} \label{simpleq}
Let Assumptions \ref{weights}, \ref{connectivity}, and
\ref{quantizedinitials}  hold.
Let $\{{x(k)}\}$ be generated by the update rule (\ref{quantupdate}).
If $k \geq nBK$, {then} all {components} of $x(k)$
are equal.
\end{proposition}

\begin{proof} Consider the nodes {whose initial
value is $U$.} There are at most $n$ of them. {As long as} not all
entries of ${x(k)}$ are equal, then every $B$ iterations, at least
one {node} must use a value strictly less than $U$ {in} an update;
\an{such a} node will have its value decreased to $U-1/Q$ or less. It
follows that after $nB$ iterations, the largest node {value will be}
at most $U-1/Q$. Repeating this argument, we {see} that at most
$nBK$ iterations are possible before {all the nodes have} the same
value.
\end{proof}

Although the above bound gives informative results for small $K$, it
becomes weaker as {$Q$ (and, therefore, $K$)} increases.
On the other hand, as $Q$ approaches infinity,
the quantized system approaches the unquantized system; the availability
of convergence time bounds for the unquantized system suggests that
similar bounds should be possible for the quantized one. Indeed,  in
the next subsection, \ao{ we adopt a notion of convergence time
parallel to our notion of convergence time for the unquantized
algorithm; as a result, we} obtain a bound on the convergence time
{which is} independent of the total number of quantization levels.

\subsection{A quantization level independent bound}

We adopt a slightly different measure of convergence for the
analysis of the quantized consensus algorithm. For any $x\in\R^n$,
we define $m(x)=\min_i x_i$ and
\[\underline{V}(x)=\sum_{i=1}^n (x_i - m(x))^2.\]
We will also use the simpler notation $m(k)$ and $\underline{V}(k)$
to denote $m(x(k))$ and $\underline{V}(x(k))$, respectively, where
it is more convenient to do so. The function $\underline{V}$ will be
our Lyapunov function for the analysis of the quantized consensus
algorithm. The reason for not using our earlier Lyapunov function,
$V$, is that for the quantized algorithm, $V$ is not guaranteed to
be monotonically nonincreasing in time. On the other hand, we have
that $V(x)\leq \underline{V}(x) \leq \alexo{4} \ao{n}V(x)$ for
any\footnote{\alexo{The first inequality follows \jnt{because}
$\sum_{i} (x_i - z)^2$ is minimized when $z$ is the mean of the
vector $x$; to establish the second inequality, observe that it
suffices to consider the case when the mean of $x$ is $0$ and
$V(x)=1$. In that case,  the largest distance between $m$ and any
$x_i$ is $2$ by the triangle inequality, so $\underline{V}(x) \leq
4n$.}}  $x \in \R^n$. As a consequence, any convergence time bounds
expressed in terms of $\underline{V}$ translate to essentially the
same bounds expressed in terms of $V$, \ao{up to a logarithmic
factor}.

Before proceeding, we {record an elementary fact which will allow us
to} relate the variance decrease $V(x)-V(y)$ to the decrease,
$\underline{V}(x)- \underline{V}(y)$, of our new Lyapunov function.
The proof involves {simple algebra,} and is therefore omitted.

\begin{lemma} \label{sumlemma} Let $u_1,\ldots,u_{n}$ and
$w_1,\ldots,w_n$ be real numbers satisfying \[ \sum_{i=1}^n u_i =
\sum_{i=1}^n w_i.
\] {Then, the expression}
\[ f(z) = \sum_{i=1}^n (u_i - z)^2 - \sum_{i=1}^n (w_i - z)^2 \]
is {a constant,}
independent of the scalar $z$.
%\[ f(z_1) = f(z_2),\quad \mbox{ for all } z_1,z_2. \]
\end{lemma}

Our next lemma places a bound on {the decrease of} the Lyapunov
function $\underline{V}(t)$ between times $kB$ and $(k+1)B-1$.

\noindent \begin{lemma} \label{quantdiff}
Let Assumptions
\ref{weights}, \ref{weakconnect}, and \ref{quantizedinitials} hold.
Let $\{x{(k)}\}$ be generated by the update rule
(\ref{quantupdate}). Suppose that the components $x_i(kB)$ of
the vector $x(kB)$ have been ordered from largest to smallest,
with ties broken arbitrarily. Then, we have
\[ \underline{V}(kB) - \underline{V}((k+1)B)
\geq \frac{\eta}{2}\, \sum_{i=1}^{n-1} (x_{i}(kB)
- x_{i+1}(kB))^2.\]
\end{lemma}

\begin{proof} For all $k$, we view Eq.\ (\ref{quantupdate}) as the
composition of two operators: \[ y(k) = A(k) x(k), \] where $A(k)$
is {a} doubly stochastic matrix, and \[ x(k+1) = \lfloor y(k)
\rfloor,
\] where the quantization $\lfloor \cdot \rfloor$ is {carried out}
componentwise.

We apply Lemma \ref{sumlemma} with the identification
$u_i=x_{{i}}(k)$, $w_i=y_{{i}}(k)$. Since multiplication by a doubly
stochastic matrix preserves the mean, the condition $\sum_i u_i =
\sum_i w_i$ is satisfied. By considering two different choices for
the scalar $z$, namely, $z_1=\bar{x}(k)=\bar{y}(k)$ and $z_2=m(k)$, we obtain
\begin{equation} \label{firststep} V(x(k)) - V(y(k)) =
\underline{V}(\ao{x(k)}) - \sum_{i=1}^n (y_i(k)-m(k))^2.
\end{equation}
Note that $x_i(k+1)-m(k)\leq y_i(k)-m(k)$. Therefore,
\begin{equation} \label{secondstep} \underline{V}(\ao{x(k)}) -
\sum_{i=1}^n (y_i(k)-m(k))^2 \leq \underline{V}(\ao{x(k)}) -
\sum_{i=1}^n (x_i(k+1)-m(k))^2.
\end{equation}
Furthermore, note that \alexo{since $x_i(k+1) \geq m(k+1) \geq m(k)$
for all $i$, we have that } $x_i(k+1)-m(k+1)\leq x_i(k+1)-m(k)$.
Therefore,
\begin{equation} \label{thirdstep} \underline{V}(\ao{x(k)}) -
\sum_{i=1}^n (x_i(k+1)-m(k))^2 \leq \underline{V}(\ao{x(k)}) -
\underline{V}(\ao{x(k+1)}).
\end{equation}
By combining
Eqs.\ (\ref{firststep}), (\ref{secondstep}), and (\ref{thirdstep}),
we obtain
\[ V(x(t)) - V(y(t)) \leq \underline{V}(\ao{x(t)}) -
\underline{V}(\ao{x(t+1)})\qquad \hbox{for all }t. \]
Summing the
preceding relations over $t=kB,\ldots, (k+1)B-1$, we further obtain
\[%\begin{equation}
\sum_{t=kB}^{(k+1)B-1} \Big(V(x(t))-V(y(t))\Big) \leq
\underline{V}(\ao{x(kB)}) - \underline{V}(\ao{x((k+1)B))}.
\]%\label{sumvariance}\end{equation}

To complete the proof, we provide a lower bound on the expression
\[
\sum_{t=kB}^{(k+1)B-1} \Big(V(x(t)) - V(y(t))\Big). \] Since
$y(t)=A(t)x(t)$ for all $t$, it follows from Lemma \ref{vl} that for
any $t$, \[ V(x(t))-V(y(t)) = \sum_{i<j} w_{ij}(t)
(x_i(t)-x_j(t))^2,\] where $w_{ij}(t)$ is the $(i,j)$-th entry of
$A(t)^T A(t)$. Using this relation and following the same line of
analysis used in the proof of Lemma \ref{vardiff} [where the
relation $\alexo{x_i(t) \geq y_{d_1}}$
%in Eq.\ (\ref{cutbound})
holds in view of the assumption that $x_i(kB)$ is a multiple of
$1/Q$ for all $k\ge 0$, cf.\ Assumption \ref{quantizedinitials}] ,
we obtain the desired result.
\end{proof}

\vspace{1pc} The next theorem contains our main result on the
convergence time of the quantized algorithm.

\noindent \begin{theorem} \label{qbound} Let Assumptions
\ref{weights}, \ref{weakconnect}, and \ref{quantizedinitials} hold.
Let $\{x{(k)}\}$ be generated by the update rule
(\ref{quantupdate}). Then, there exists an absolute constant $c$ such that
we have
$$\underline{V}(k) \leq \epsilon \underline{V}(0)\qquad \hbox{for all }
k\ge c\, (n^2/\eta) B
\log(1/\epsilon).$$

\end{theorem}

\begin{proof}
Let us assume that $\underline{V}(kB)>0$. From Lemma
\ref{quantdiff}, we have
\[ \underline{V}(kB)- \underline{V}((k+1)B) \geq \frac{\eta}{2}
\sum_{i=1}^{n-1} (x_{i}(kB) - x_{i+1}(kB))^2,  \] where {the
components} $x_i(kB)$ are ordered from largest to smallest. Since
$\underline V(kB) = \sum_{i=1}^n (x_i(kB) - x_n(kB))^2$, we have \[
\frac{\underline{V}(kB)-\underline{V}((k+1)B)}{\underline{V}(kB)}
\geq \frac{\eta}{2} \frac{\sum_{i=1}^{n-1} (x_{i}(kB) -
x_{i+1}(kB))^2}{\sum_{i=1}^n (x_i(kB)-x_n(kB))^2}.
\] {Let} $y_i = x_i(kB)-x_n(kB)$. Clearly, $y_i \geq 0$ for all $i$,
and $y_n=0$. Moreover, the monotonicity of $x_i(kB)$ implies the
monotonicity of $y_i$: \[ y_1 \geq y_2 \geq \cdots \geq y_n=0.
\] Thus,
\[ \frac{\underline{V}(kB)-\underline{V}((k+1)B)}{\underline{V}(kB)}
\geq \frac{\eta}{2} \min_{{y_1 \geq y_2\ge \cdots\ge y_n\atop
y_n=0}}\frac{\sum_{i=1}^{n-1} (y_i - y_{i+1})^2}{\sum_{i=1}^n
y_i^2}.
\] \alexo{Next, we simply repeat the steps of Lemma \ref{lboundvar}. We can
assume without loss of generality that $\sum_{i=1}^n y_i^2=1$.
\an{Define $z_i=y_{i}-y_{i+1}$ for $i=1,\ldots,n-1$ and $z_n=0$. We}
have that $z_i$ are all nonnegative
and $\sum_{i} z_i = y_1 - y_n \geq 1/\sqrt{n}$.} \an{Therefore,}
\[ \frac{\eta}{2} \min_{{y_1 \geq y_2\geq \cdots \geq y_n\atop
\sum_i y_i^2=1}} \sum_{i=1}^{n-1} (y_i - y_{i+1})^2 \geq
\frac{\eta}{2} \min_{z_i \geq 0, \sum_i z_i \geq 1/\sqrt{n}}
\sum_{i=1}^n z_i^2.
\]
The minimization problem on the right-hand side has an optimal value of at
least $1/n^2$, and the desired result follows.
\end{proof}

\subsection{Extensions and modifications}

In this subsection, we  comment briefly on some corollaries
of Theorem \ref{qbound}.

First, we note that the results of Section \ref{matrixpicking}
immediately carry over to the quantized case. Indeed, in Section
\ref{matrixpicking}, we showed how to pick the weights $a_{ij}(k)$
in a decentralized manner, based  only on local information, so that
Assumptions 1 and \ref{weakconnect} are satisfied, with $\eta \geq
1/3$. When using a quantized version of the balancing algorithm,
%of Section \ref{matrixpicking},
we once again \alexo{ manage to remove the factor of $1/\eta$ from
our upper bound.}

\begin{proposition} \label{savingnq}
For the quantized version of the balancing algorithm,
%of Section \ref{matrixpicking},
and under the same assumptions as in Theorem \ref{savingn}, if
$k\geq c\, n^2B \log (1/\epsilon))$, then $\underline{V}(k ) \leq
\epsilon \underline{V}(0)$, where $c$ is an absolute
constant.\end{proposition}

Second, we note that Theorem \ref{qbound} {can be used} to obtain a
bound on the time until {the values of all nodes are equal. Indeed,
we observe that in the presence of quantization, once the condition
$\underline{V}(k) <1/Q^2$ is satisfied, all components of $x(k)$
must be equal.}

\begin{proposition} \label{qboundequal}
Consider the quantized algorithm (\ref{quantupdate}), and assume
that Assumptions \ref{weights}, \ref{weakconnect}, and
\ref{quantizedinitials} hold. If $k\geq  c(n^2/\eta) B \big[\log Q +
\log \underline{V}(0)\big]$, then all components of $x(k)$ are
equal, where $c$ is an absolute constant. \end{proposition}

\subsection{Tightness}

{We now} show that the \alexo{quantization-level independent} bound
in Theorem \ref{qbound} is tight, \alexo{even when the weaker}
Assumption \ref{weakconnect} is replaced with \alexo{the stronger}
Assumption \ref{connectivity}.

\begin{proposition} \label{quanttight}
There \jnt{exist} absolute constant\alexo{s} $c$ \alexo{ and $n_0$
}with the following property. For {any} nonnegative integer $B$,
$\eta < 1/2$, $\epsilon < 1$, and \alexo{ and $n \geq n_0$ }, there
exist a sequence of weight matrices $A(k)$ satisfying Assumptions
\ref{weights} and \ref{connectivity}, and an initial value $x(0)$
satisfying Assumption \ref{quantizedinitials}, \alexo{ and a number
quantization levels $Q(n)$ (depending on $n$)} such that under the
dynamics of Eq.\ (\ref{quantupdate}), {if}
$\underline{V}(k)/\underline{V}(0) \leq \epsilon$, {then}
\[ k\geq c\, \frac{n^2}{\eta}\, B \log \frac{1}{\epsilon}.
\]
\end{proposition}
\begin{proof} We have demonstrated in \alexo{Proposition}
\ref{unquanttight} a similar result for the unquantized algorithm.
Namely, we have shown that for {$n$ \ao{large enough} and for any}
$B$, $\eta<1/2$, and $\epsilon<1$, there exists a weight sequence
$a_{ij}(k)$ and an initial {vector} $x(0)$ such that the first time
when $V(t) \leq \epsilon V(0)$ occurs after $\Omega((n^2/\eta) B
\log (1/\epsilon))$ steps. Let $T^*$ be this first time.

{Consider the quantized algorithm under} the exact same sequence
$a_{ij}(k)$, {initialized at $\lfloor x(0) \rfloor$.}  Let
$\hat{x}_i(t)$ refer to the value of node $i$ at time $t$ in the
quantized algorithm {under} this scenario, {as opposed to} $x_i(t)$
which {denotes the} value in the unquantized algorithm. Since
quantization can {only} decrease a nodes value by at most $1/Q$ at
each iteration, it is easy to show, by induction, that
\[ x_i(t) \geq \hat{x}_i(t) \geq x_i(t) - t/Q \] We can pick $Q$
large enough so that, {for $t < T^*$,} the vector $\hat{x}(t)$ is
{as close as desired to} $x(t)$.

Therefore, for $t < T^*$ and for large enough $Q$,
$\underline{V}(\hat{x}(t))/\underline{V}(\hat{x}(0))$ will be
arbitrarily close to $\underline{V}(x(t))/\underline{V}(x(0))$.
{From the proof of \alexo{Proposition} \ref{unquanttight}, we see
that} $x(t)$ is always a {scalar} multiple of $x(0)$. Since
$\underline{V}(x)/V(x)$ is invariant under multiplication by a
constant, it follows that $\underline{V}(x(t))/\underline{V}(x(0)) =
V(x(t))/V(x(0))$. Since this last quantity is above $\epsilon$ for
$t<T^*$, it follows that provided $Q$ is large enough,
$\underline{V}(\hat{x}(t))/\underline{V}(\hat{x}(0))$ is also above
$\epsilon$ for $t<T^*$. This proves the proposition. \end{proof}

\subsection{Quantization error}

Despite favorable convergence properties of our quantized
{averaging} algorithm (\ref{quantupdate}), the update rule does not
preserve the average of the values at each iteration. Therefore, the
{common limit of the sequences $x_i(k)$, denoted by $x_f$,}
need not be equal to the exact average of the initial values. We
next provide an upper bound on the error between {$x_f$} and the
initial average, as a function of the number of quantization levels.

\begin{proposition} \label{qerror}
There is an absolute constant $c$ such that
for the common limit $x_f$ of the
values $x_i(k)$ generated by the
quantized algorithm (\ref{quantupdate}), we have
\[ \left|x_f - \frac{1}{n} \sum_{i=1}^n x_i(0)\right| \le {c\over Q}
\ {n^2\over\eta}\, B\log (Qn(U-L)).\]
\end{proposition}

\begin{proof} By \alexo{Proposition} \ref{qboundequal}, after $O\Big((n^2/\eta) B \log (Q
\underline{V}(x(0)))\Big)$ iterations, all nodes will have the same
value. Since $\underline{V}(x(0)))\le n(U-L)^2$ and the average
decreases by at most $1/Q$ at each iteration, the result follows.
\end{proof}

\vspace{1pc}

%This theorem quantifies the relationship between network size and the number of quantization levels per unit $Q$. If we are interested in how the performance of an averaging algorithm scales with the network size $n$, then this theorem implies that the number of bits, which is proportional to $\log Q$, needs to grow only as $O(\log n)$ to keep the precision of the quantized algorithm the same.

%\noindent {\bf Remark}:
Let us assume that the parameters %of the problem
$B$, $\eta$, and $U-L$ are fixed. \alexo{Proposition} \ref{qerror}
implies that as $n$ increases, the number of bits {used for each
communication,} which is proportional to $\log Q$, needs to grow
only as $O(\log n)$ to make the error negligible. Furthermore, this
is true even if the {parameters} $B$, ${1}/{\eta}$, and $U-L$ grow
polynomially in $n$.

\begin{center}
\begin{figure}
\hspace{2cm}
\includegraphics*[width=12cm]{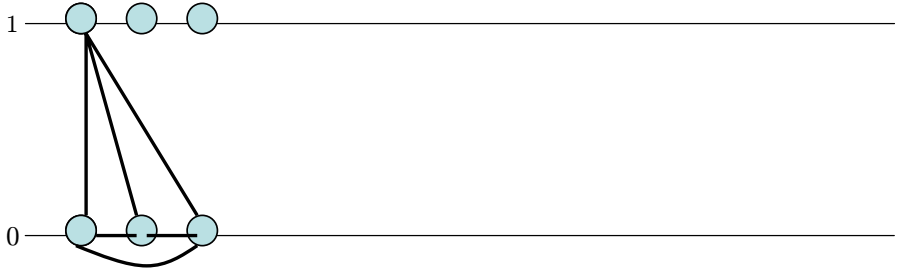}
\caption{\label{ic} Initial configuration. Each node takes the
average value of its neighbors. }
\end{figure} \end{center}

{For a converse, it can} be seen that ${\Omega}(\log n)$ bits are
needed. Indeed, consider  $n$ nodes, with $n/2$ {nodes initialized
at} $0$, and $n/2$ {nodes initialized at} $1$. {Suppose} that ${Q} <
n/2$; we connect the nodes by {forming} a complete subgraph over all
the nodes with value $0$ and exactly {one} node with value $1$; see
Figure \ref{ic} for an example with $n=6$. Then, each node {forms
the} average {of} its neighbors. This brings one of the nodes with
{an initial value of} $1$ down to $0$, without raising the value of
any {other} nodes. We can repeat this {process,} to bring all of the
nodes with {an initial value of} $1$ down to $0$. Since the true
average is $1/2$, the final result is $1/2$ away from the true
average. {Note now that $Q$ can} grow linearly with $n$, and still
satisfy the inequality ${Q}<n/2$. {Thus,} the number of bits can
grow as $\Omega(\log n)$, and yet, independent of $n$, the error
remains $1/2$.

\section{Conclusions} \label{conclusions}

We studied distributed algorithms for the {averaging} problem over
networks with time-varying topology, {with a focus on} tight bounds
on the convergence time  of a general class of {averaging}
algorithms. We first considered {algorithms for the case where}
agents can exchange and store continuous values, {and established
tight convergence time bounds.} We next studied averaging algorithms
under the additional constraint that agents can {only} store and
send quantized values. We showed that these algorithms guarantee
convergence of the agents values to consensus within some error from
the average of the initial values. We provided a bound on the error
that highlights the dependence on the number of quantization levels.

Our paper is a contribution to the growing literature on distributed
control of multi-agent systems. Quantization effects are an integral
part of such systems but, with the exception of a few recent
studies, have not attracted much attention in the vast literature on
this subject. In this paper, we studied a quantization scheme that
guarantees consensus at the expense of some error from the initial
average value. We used this scheme to study the effects of the
number of quantization levels on the convergence time of the
algorithm and the distance from the true average.

The framework provided in this paper motivates a number of further
research directions:
\begin{itemize}
\item[(a)] The algorithms studied in this paper assume that there is no delay in receiving the values of the other agents,
 which is a restrictive assumption in network settings. Understanding the convergence of averaging algorithms and implications
  of quantization in the presence of delays is an important topic for future research.
\item[(b)] We studied a quantization scheme with favorable convergence properties, that is, rounding down to the nearest quantization level.
    Investigation of other quantization schemes and their impact on convergence time and error is left for future work.
\item[(c)] The quantization algorithm we adopted implicitly assumes that the agents can carry out
 computations with continuous values, but can store and transmit only quantized values.
 Another interesting area for future work is to incorporate the additional constraint of finite precision computations into the quantization scheme.
 \item[(d)] \alexo{Although our bounds are tight in the worst case
 over all graphs, they are not guaranteed to perform better on
 well-connected graphs as compared to sparse graphs with many
 potential bottlenecks. An interesting question is whether it \jnt{is} be possible to pick averaging
 algorithms that learn the graph and make optimal use of its information diffusion properties. }
\end{itemize}

\newpage

\begin{small}

\bibliographystyle{amsplain}
\bibliography{distributed_alex}

\providecommand{\bysame}{\leavevmode\hbox to3em{\hrulefill}\thinspace}
\providecommand{\MR}{\relax\ifhmode\unskip\space\fi MR }
\providecommand{\MRhref}[2]{%
  \href{http://www.ams.org/mathscinet-getitem?mr=#1}{#2}
}
\providecommand{\href}[2]{#2}
\begin{thebibliography}{10}

\bibitem{BT89}
D.P. Bertsekas and J.N. Tsitsiklis, \emph{Parallel and distributed computation:
  Numerical methods}, Prentice Hall, 1989.

\bibitem{BFT05}
P.A. Bliman and G.~Ferrari-Trecate, \emph{Average consensus problems in
  networks of agents with delayed communications}, Proceedings of the Joint
  44th IEEE Conference on Decision and Control and European Control Conference,
  2005.

\bibitem{BHOT05}
V.D. Blondel, J.M. Hendrickx, A.~Olshevsky, and J.N. Tsitsiklis,
  \emph{Convergence in multiagent coordination, consensus, and flocking},
  Proceedings of the Joint 44th IEEE Conference on Decision and Control and
  European Control Conference, 2005.

\bibitem{boyd}
S.~Boyd, A.~Ghosh, B.~Prabhakar, and D.~Shah, \emph{Gossip algorithms: Design,
  analysis and applications}, Proceedings of IEEE INFOCOM, 2005.

\bibitem{CFFTZ07}
R.~Carli, F.~Fagnani, P.~Frasca, T.~Taylor, and S.~Zampieri, \emph{Average
  consensus on networks with transmission noise or quantization}, Proceedings
  of European Control Conference, 2007.

\bibitem{CFSZ05}
R.~Carli, F.~Fagnani, A.~Speranzon, and S.~Zampieri, \emph{Communication
  constraints in the state agreement problem}, Preprint, 2005.

\bibitem{C06}
J.~Cortes, \emph{Analysis and design of distributed algorithms for
  chi-consensus}, Proceedings of the 45th IEEE Conference on Decision and
  Control, 2006.

\bibitem{C89}
G.~Cybenko, \emph{Dynamic load balancing for distributed memory
  multiprocessors}, Journal of Parallel and Distributed Computing \textbf{7}
  (1989), no.~2, 279--301.

\bibitem{G06}
R.M. Gray, \emph{Toeplitz and circulant matrices: A review}, Foundations and
  Trends in Communications and Information Theory \textbf{2} (2006), no.~3,
  155--239.

\bibitem{JLM03}
A.~Jadbabaie, J.~Lin, and A.S. Morse, \emph{Coordination of groups of mobile
  autonomous agents using nearest neighbor rules}, IEEE Transactions on
  Automatic Control \textbf{48} (2003), no.~3, 988--1001.

\bibitem{KBS06}
A.~Kashyap, T.~Basar, and R.~Srikant, \emph{Quantized consensus}, Proceedings
  of the 45th IEEE Conference on Decision and Control, 2006.

\bibitem{LO81}
H.J. Landau and A.M. Odlyzko, \emph{Bounds for the eigenvalues of certain
  stochastic matrices}, Linear Algebra and its Applications \textbf{38} (1981),
  5--15.

\bibitem{LR06}
Q.~Li and D.~Rus, \emph{Global clock synchronization in sensor networks}, IEEE
  Transactions on Computers \textbf{55} (2006), no.~2, 214--224.

\bibitem{M05}
L.~Moreau, \emph{Stability of multiagent systems with time-dependent
  communication links}, IEEE Transactions on Automatic Control \textbf{50}
  (2005), no.~2, 169--182.

\bibitem{OSM04}
R.~Olfati-Saber and R.M. Murray, \emph{Consensus problems in networks of agents
  with switching topology and time-delays}, IEEE Transactions on Automatic
  Control \textbf{49} (2004), no.~9, 1520--1533.

\bibitem{OT06}
A.~Olshevsky and J.N. Tsitsiklis, \emph{Convergence rates in distributed
  consensus averaging}, Proceedings of the 45th IEEE Conference on Decision and
  Control, 2006.

\bibitem{RB05}
W.~Ren and R.W. Beard, \emph{Consensus seeking in multi-agent systems under
  dynamically changing interaction topologies}, IEEE Transactions on Automatic
  Control \textbf{50} (2005), no.~5, 655--661.

\bibitem{T84}
J.N. Tsitsiklis, \emph{Problems in decentralized decision making and
  computation}, Ph.D. thesis, Dept. of {E}lectrical {E}ngineering and
  {C}omputer {S}cience, {MIT}, 1984.

\bibitem{TBA86}
J.N. Tsitsiklis, D.P. Bertsekas, and M.~Athans, \emph{Distributed asynchronous
  deterministic and stochastic gradient optimization algorithms}, IEEE
  Transactions on Automatic Control \textbf{31} (1986), no.~9, 803--812.

\bibitem{VCBJCS95}
T.~Vicsek, A.~Czirok, E.~Ben-Jacob, I.~Cohen, and O.~Schochet, \emph{Novel type
  of phase transitions in a system of self-driven particles}, Physical Review
  Letters \textbf{75} (1995), no.~6, 1226--1229.

\bibitem{XB04}
L.~Xiao and S.~Boyd, \emph{Fast linear iterations for distributed averaging},
  Systems and Control Letters \textbf{53} (2004), 65--78.

\bibitem{XBK07}
L.~Xiao, S.~Boyd, and S.-J. Kim, \emph{Distributed average consensus with
  least-mean-square deviation}, Journal of Parallel and Distributed Computing
  \textbf{67} (2007), 33--46.

\end{thebibliography}

\end{small}
\end{document}